\documentclass[a4paper,10pt]{amsart}

\usepackage{amssymb}
\usepackage{amsfonts}
\usepackage{amsmath}
\usepackage{graphicx}%
\usepackage{lastpage}
\usepackage[utf8]{inputenc}
\usepackage[legalpaper,bookmarks=true,colorlinks=true,linkcolor=blue,citecolor=blue]{hyperref}
\usepackage{fancyhdr}
\usepackage{color}
\usepackage[mathlines]{lineno}
\usepackage{multirow}
\usepackage{caption}
\usepackage{subcaption}
\usepackage{lscape}
\usepackage{epsfig}
\usepackage{natbib}
\usepackage{float}
\usepackage{slashbox} % Defines \backslashbox

\newtheorem{theorem}{Theorem}
\theoremstyle{plain}

\newtheorem{corollary}{Corollary}

\newtheorem{definition}{Definition}

\numberwithin{equation}{section}
\begin{document}

\title{
Entropies and their Asymptotic Theory in the  Discrete case }
\author{$^{(1)}$ 
BA Amadou Diadi\'e}
\email{ba.amadou-diadie@ugb.edu.sn}
\author{$^{(1,2,4)}$ LO Gane Samb}
\email{gane-samb.lo@ugb.edu.sn}

\begin{abstract}
We present some new nonparametric estimators of entropies and we establish almost sure consistency and central limit Theorems for some of the most important entropies in the discrete case. Our theorical results are validated by simulations. 
\end{abstract}

\maketitle
\section{Introduction}
\subsection{Motivation}

Consider an outcome $A$ of a random experiment on a probability space $(\Omega,\mathcal{A},\mathbb{P})$. The \textit{information amount
} or \textit{content} of the outcome $A$ is (see \cite{carter}) $$\mathcal{I}(A)=\log_2 \frac{1}{\mathbb{P}(A)},$$where $\log_2$ is the logarithm base $2$.\\

\noindent Let $X$ be a discrete random variable 
defined on a probability space $(\Omega,\mathcal{A},\mathbb{P})$ and $\{c_1,c_2,\cdots,c_r\}$ the set of all possible values of $X$. \\

\noindent The probability
distribution $\textbf{p}=(p_j)_{(j=1,\cdots,r)}$ of the events $(X=c_j)
%,\{X=c_2\},\cdots,\{X=c_r\}
$, coupled with the \textit{information
amount} of every event $\mathcal{I}(X=c_j)_{(j=1,\cdots,r)}$, forms a random variable whose
expected value is the \textit{average amount of information}, or
\textit{entropy} (more specifically, \textit{Shannon entropy}) generated by this distribution.
\begin{definition}
Let $X$ be a discrete random variable defined on the probability space $(\Omega,\mathcal{A},\mathbb{P})$ and taking values in the finite countable space $\mathcal{X}=\{c_1,c_2,\cdots,c_r\}\,(r\geq 2)$ with probabilities mass function (p.m.f.) %$\textbf{p}=(p_c)_{c\in \mathcal{X}}$ 
$(\textbf{p})_{p_j,j=1,\cdots,r}$, that is,  $ p_{j}=\mathbb{P}(X= c_j)\ \ \forall j\in J= \{1,\cdots,r\}.$\\

\noindent  The \textit{Shannon entropy} of the random variable $X$ is given by 
\begin{equation}
\label{extshan}
\mathcal{E}_{Sh}(X
)=\sum_{j=1}^r  p_j
%\mathcal{I}(\{X=c_j\})
\log_2\frac{1}{p_j}=\mathbb{E}\left( \log_2 (\textbf{p})\right).
\end{equation}
\end{definition}

\noindent Entropy is usually measured in \textit{bits} (\textbf{b}inary \textbf{i}nformation uni\textbf{t}) (if $\log_2$), \textit{nats} (if $\ln )$, or \textit{hartley}( if $\log_{10})$, depending on the base
of the logarithm which is used to define it.\\

\noindent For ease of computations and notation convenience, we use the natural logarithm ($\ln$) since logarithms of varying bases are related by a constant.\\

\bigskip \noindent In the sequel, we consider the entropy of the discrete random variable $X$ as a function of discretes probabilities  $\textbf{p}=(p_j)_{j\in J}$.

\subsection{Generalizations of Shannon entropy}
\noindent Inspired by the study of \textit{$\alpha$-deformed algebras} and
special functions, various generalizations have been investigated. \\

\noindent Most notably, \cite{ren2} proposed a one parameter family of
entropies extending \textit{Shannon entropy}.\\

\noindent (b) 
The $\alpha-$\textit{R\'enyi entropy} of the random variable $X$ is defined by% (\cite{ren1})
\begin{equation}
\label{reny} \mathcal{E}_{R,\alpha} (\textbf{p})=\frac{1}{1-\alpha} \ln \left( \sum_{j=1}^rp_{j}^\alpha\right).
\end{equation}

\noindent with $\alpha\in (0,1)\cup (1,+\infty)$, which, in particular, reduces to the \textit{Shannon entropy} in the limit $\alpha\rightarrow 1$.\\

\noindent (c) 
Also, the $\alpha-$\textit{Tsallis entropy}  of the random variable $X$ defined by (see \cite{tsal1}) :

\begin{equation}\label{stal}
\mathcal{E}_{T,\alpha} (\textbf{p})=\frac{1}{1-\alpha}\left( \sum_{j=1}^rp_{j}^\alpha-1\right),\ \ \alpha\in (0,1)\cup (1,+\infty)
\end{equation}
 has generated a large burst of research activities.\\
%fich=lec1-entropy

\bigskip \noindent Let cite a few other examples of entropies. \\

\bigskip \noindent (d) The $\alpha-$\textit{Landsberg-Vedral entropy} also called \textit{normalized Shannon entropy} of the random variable $X$ is defined by (see \cite{land}) :

\begin{equation}\label{lventrop}
\mathcal{E}_{L.V,\alpha} (\textbf{p})=\frac{1}{1-\alpha}\left( 1-\frac{1}{ \sum_{j=1}^r p_{j} ^\alpha}\right)=\frac{\mathcal{E}_{T,\alpha}(\textbf{p})}{ \sum_{j=1}^r p_{j} ^\alpha},\ \ \alpha\in (0,1)\cup (1,+\infty).
\end{equation}

\bigskip 
\noindent (e) The $\alpha-$\textit{Abe entropy}  of the random variable $X$ is defined by (see \cite{abe}) :$$\mathcal{E}_{Ab,\alpha}( \textbf{p })=-\frac{1}{\alpha-\alpha^{-1}} \sum_{j=1}^r(p_{j} ^\alpha- p_{j} ^{\alpha^{-1}}),\ \ \alpha\in (0,1)\cup (1,+\infty). $$

\bigskip 
\noindent (f) The $\kappa$-entropy of the random variable $X$
is defined by the following expression (see \cite{kana}) :
$$\mathcal{E}_{\kappa}( \textbf{p} )=\frac{1}{2\kappa } \sum_{j=1}^r(p_{j} ^{1-\kappa}- p_{j} ^{1+\kappa})
,\ \ \kappa\in (0,1).$$

\bigskip \noindent (g) The \textit{Varma's entropy} of order $\alpha$ and type $\beta$ of the random variable $X$ is defined by 

\begin{equation}
\mathcal{E}_{V,\alpha,\beta}(\textbf{p})=\frac{1}{\beta-\alpha}\ln \left( \sum_{j=1}^{r} p_{j} ^{\alpha+\beta-1}\right),\ \ \text{for}\ \ \beta-1<\alpha<\beta,\ \ \beta\geq 1
.\end{equation}
%\bigskip \noindent (h) The family of generalized entropies $\alpha,\beta-$\textit{Sharma-Mital entropy} of the random variable $X$ is given by (see \cite{shar}) :
%\begin{equation}\label{sharm}
%\mathcal{E}_{S.M,\alpha,\beta} (\textbf{p})= \frac{1}{1-\beta}\left[ \left( \sum_{j=1}^{r} p_{j} ^{\alpha} \right)^{\left( \frac{1-\beta}{1-\alpha} \right)}-1 \right]
%\end{equation}
%where $\alpha,\beta\in (0,1)\cup(1,+\infty)$ and $\alpha\neq \beta$.\\
%
%\noindent This bi-parametric family of entropies tends in limit cases to Renyi entropies (for $\beta\rightarrow 1$), Tsallis entropies (for $\beta \rightarrow \alpha$), Shannon entropy (for both $\alpha,\beta \rightarrow 1$), and the Landsberg-vedral entropy (for $\beta=2-\alpha$). \\

%
%In the limit $\beta\rightarrow 1$, Sharma-Mittal entropy becomes R\'enyi entropy,  
%while for $\beta\rightarrow \alpha$, it is Tsallis entropy.
%In the limiting case when both parameters approach $1$, the Shannon entropy is recovered.\\

\bigskip  \noindent Interestingly, the \textit{Landsberg-Vedral} and \textit{$\kappa$} entropies reduce  to the \textit{Shannon entropy} in the limit $\alpha\rightarrow 1$ and $\kappa\rightarrow 0$ respectively.\\

\noindent From this small sample of entropies, we may give the following remarks :\\

\noindent (a) For most entropies, we may have computation problems. So without loss of generality, suppose
\begin{equation}
\ \ \ p_{j}>0, \ \ \forall j\in J=\{1,\cdots,r\}
 \ \ \ \ \ (\textbf{BD}).
\label{BD}
\end{equation}

\bigskip \noindent If Assumption (\ref{BD}) holds, we do not have to worry about summation problems, especially for entropies cited below in the computations arising in estimation theories. This explain why Assumption \ref{BD} is systematically used in a great number of works in that topics, for example, in \cite{amad_glo}, \cite{singh}, \cite{kris}, \cite{hall}, to cite a few.\\

\noindent (b) 
%Most of entropies are build on the \textit{power sum} of order $\alpha\in (0,1)\cup (1,+\infty)$ of $\textbf{p}$ over $\{c_j,j\in J\}$  
%\begin{equation}
%\mathcal{S}_{\alpha}(\textbf{p})=\sum_{j=1}^r p_{j} ^\alpha.
%\end{equation}
 The \textit{power sum} of order $\alpha\in (0,1)\cup (1,+\infty)$ of the distribution \textbf{p} over $\{c_j,j\in J\}$ is
\begin{equation}\label{ialpha}\mathcal{S}_\alpha(\textbf{p})=\sum_{j\in J}p_j^\alpha,
 \end{equation} and, is related to  Reyni, Tsallis, Landsberg-Vedral, Abel, $\kappa$, and Varma entropies via
\begin{eqnarray*}
&& \mathcal{E}_{R,\alpha} (\textbf{p})=\frac{1}{1-\alpha} \ln \left( \mathcal{S}_\alpha(\textbf{p})\right),\ \ \ \ 
\mathcal{E}_{T,\alpha} (\textbf{p})=\frac{1}{1-\alpha}\left( \mathcal{S}_\alpha(\textbf{p})-1\right) \\
&& \mathcal{E}_{L.V,\alpha} (\textbf{p})=\frac{1}{1-\alpha}\left( 1-\frac{1}{  \mathcal{S}_\alpha(\textbf{p})}\right),\ \ \
\mathcal{E}_{Ab,\alpha}( \textbf{p })=-\frac{1}{\alpha-\alpha^{-1}} \left(  \mathcal{S}_\alpha(\textbf{p})+ \mathcal{S}_{\alpha^{-1}}(\textbf{p})\right), \\ 
&&\mathcal{E}_{\kappa}( \textbf{p} )=\frac{1}{2\kappa } (\mathcal{S}_{1-\kappa}(\textbf{p})-\mathcal{S}_{1+\kappa}(\textbf{p})),\ \  \text{and}\ \  \mathcal{E}_{V,\alpha,\beta}(\textbf{p})\frac{1}{\beta-\alpha}\ln \left(\mathcal{S}_{\alpha+\beta-1}(\textbf{p})\right).
\end{eqnarray*} 
Hence establishing asymptotic limits of estimators of these ones is equivalent to establishing asymptotic limits of $
\mathcal{S}_\alpha(\widehat{\textbf{p}}_n)$.

\subsection{Bibliography and applications}
\noindent Although we are focusing on the aforementioned entropies in this paper, it is worth mentioning that there exist quite a few number of them.\\

\bigskip \noindent Let us cite for example the ones named after : Fuzzy Entropy (see \cite{luca}, \cite{bhand}, \cite{kosk}, \cite{pal}, \cite{yager}), Havrda-Charv\'at entropy (see \cite{havr}), Generalized Entropy also called $f-$divergence (see \cite{liese}, \cite{bales})% conditional R\'enyi entropy
, Frank-Daffertshofer entropy (see \cite{frank}), Kapur measure (see \cite{kapur}), Hartley entropy, min entropy and max entropy (see \cite{dodis}), collision entropy etc.\\

\noindent Recently, there have been made several successful attempts in order to categorize the various entropy classes and
their properties: \cite{hanel1}, \cite{hanel2} classified the entropies according to their asymptotic scaling. \cite{temp} studied the Generalized entropies according to group properties. \cite{biro} derived a new class of
entropies from its interaction with heat reservoir. \cite{ilic} classified the pseudo-additive entropies by
generalization of Khinchin axioms.\\

\noindent Before coming back to our entropies estimation of interest, we want to highlight some important applications of them.\\

\noindent Indeed, entropy has proven to be useful in applications. Let us cite some of them :\\

\noindent (a) The entropy concept was born initially in thermodynamics by \cite{claus} to measure the ratio of transferred heat through a reversible process in an isolated
system and to measure of
uncertainty about the system that remains after observing its macroscopic
properties (pressure, temperature or volume). Since then, entropy has been of great theoretical and applied interest.\\

\noindent b) In finance, \cite{phil} were the first two authors who applied the
concept of entropy to portfolio selection. It has been used as a risk measure for stock, for portfolio returns, for portfolio diversifications (see \cite{miha}), it has been applied as measure of investment risk in the
discrete case (see \cite{nawr}), as well as a measure of dependence in return time series (see \cite{maa}). \\

\noindent (c) While a significant number of other entropies have since been
introduced, \textit{R\'enyi entropy}  is especially important because it is a well known one parameter generalization of Shannon entropy. It is often used as a bound on Shannon entropy (see \cite{mokk}, 
\cite{nem}, \cite{har}, and it replaces Shannon entropy as a measure of randomness (see \cite{csis}, \cite{massey}, \cite{arik},etc). It generalizes also Hartley, collision, and min-entropy. It has
successfully been used in a number of different fields, such as statistical physics, quantum mechanics,
communication theory and data processing (see \cite{jizb}, \cite{csis}), in the context of channel coding (see \cite{arim}), secure communication (see \cite{cach}),  and \cite{haya}),  multifractal analysis (see \cite{jizb}). \\

\noindent In the context of fractal dimension estimation, \textit{R\'enyi entropy} forms the basis of the concept of generalized dimensions. It intervene as well in ecology and statistics as index of diversity.\\

\noindent R\'enyi entropy is also of interest in its own right, with diverse applications in unsupervised learning (see \cite{xu}, \cite{jens}, source adaptation (see \cite{mans}, image registration (see \cite{ma}, see \cite{neem}, 
 and password guessability (see \cite{arik},
\cite{pfi}, \cite{hana} among others. In particular, the R\'enyi entropy of order $2$ measures the quality
of random number generators (see \cite{knu}) , \cite{oor}, determines the number of unbiased bits that can be
extracted from a physical source of randomness see \cite{impa}, \cite{ben}, helps test graph expansion \cite{gold} and
closeness of distributions \cite{bat}, and characterizes the number of reads needed to reconstruct
a DNA sequence \cite{mota}.\\

\noindent  The \textit{R\'enyi entropy} is important in ecology and statistics as index of diversity. It  is also important in quantum information, where it can be used as a measure of entanglement. \\

\noindent e) \textit{Varma's entropy} plays a vital role as a measure of complexity and uncertainty in different
areas such as physics, electronics and engineering to describe many chaotic systems\\

\noindent f)  In the context of multi-dimensional harmonic oscillator
systems, the \textit{Sharma–Mittal
entropy} has previously been studied (see \cite{Uzengi}).

\subsection{Previous work}
The estimation of entropies have become growingly important for their wide applications in the fields of neural science and information theory, etc.\\

\noindent For example Shannon entropy estimation has several applications, including measuring genetic diversity (see \cite{shen}, quantifying neural activity (see \cite{panin}), see \cite{nem}, network anomaly detection \cite{lall}, and others.\\

\noindent Most texts on entropy estimation deal with Shannon entropy estimation and use the plug-in method. \\

\noindent \cite{zing} showed that, if $\{p_j,\,j\geq 1\}$ is non uniform distribution satisfying $\mathbb{E}(\log P_X)^2<\infty,$ and   if there exists an integer valued function $J(n)$ such that, $J(n)\rightarrow+\infty,\ \ J(n)=o(\sqrt{n})$
and
$\sqrt{n}\sum_{j\geq J(n)}p_j\log p_j\rightarrow 0,$ as $n\rightarrow \infty$, then

$$\sqrt{n}(\mathcal{E}_{Sh}(\widehat{\textbf{p}}_n)-\mathcal{E}_{Sh}(\textbf{p})\rightsquigarrow \mathcal{N}(0,\sigma_{Sh}^2(\textbf{p}))\ \ \text{as}\ \ n\rightarrow+\infty$$

\noindent where $\sigma_{Sh}^2(\textbf{p})=Var(-\log \mathbb{P}_X)>0$.\\

%
%
%fich IT-12-0260R1\\

%\noindent \cite{zhao}
%Approximated the empirical \textit{Shannon entropy} from
%estimates of R\'enyi $/$ Tsallis entropies by approximating the fonction $s\mapsto s\log s$ by a linear combination of two functions of the form $s\mapsto s^r,\,r\in(0,2]$ for $s$ less that an upper bound.

\noindent \cite{zhang} proposed a non parametric estimator of Shannon's entropy on a countable alphabet 
$$\mathcal{E
}^{(Z)}(\widehat{\textbf{p}}_n)=\sum_{\ell=1}^{n-1}\frac{1}{\ell}\left\lbrace \frac{n^{\ell+1}[n-(\ell+1)]!}{n!}\sum_j\left[\widehat{p}_n^{c_j} \prod_{i=0}^{\ell-1}\left(1-\widehat{p}_n^{c_j}-\frac{i}{n} \right)\right] \right\rbrace$$
and established that
$$\mathbb{E}\left(\mathcal{E
}^{(Z)}(\widehat{\textbf{p}}_n)
\right)-\mathcal{E}(\textbf{p})=O\left(\frac{(1-p_0)^n}{n}\right)$$ where $p_0=\min_{j\in J}\{p_j\}$.\\

\noindent Later on \cite{miller2}, \cite{basha}, and \cite{harris} established that
\begin{eqnarray}
\label{biais}&&\mathbb{E}\left( \mathcal{E}_{Sh}(\widehat{\textbf{p}}_n)-\mathcal{E}_{Sh}(\textbf{p})\right)=-\frac{r-1}{2n}+\frac{1}{12n^2}\left(1-\sum_{j=1}^r\frac{1}{p_k} \right)+O(n^{-3})\\
\label{variance} &&Var(\mathcal{E}_{Sh}(\widehat{\textbf{p}}_n))= 
\frac{1}{n}\left(\sum_{j=1}^rp_j(\ln p_j)^2-(\mathcal{E}_{Sh}(\textbf{p}))^2\right)+\frac{r-1}{2n^2}+O(n^{-3})
\end{eqnarray}
%The plug-in estimator for $\mathcal{E}(\textbf{p})$, given by 
%$$\mathcal{E}(\widehat{\textbf{p}}_n) = $$ 

\bigskip
% fichier Zhang-uncc-06\\
 
 \bigskip \noindent
\cite{anto1} proved that 
$$\mathbb{E}\left(\mathcal{E}(\widehat{\textbf{p}}_n)-\mathcal{E}(\textbf{p})\right) \sim n^{-(\lambda-1)/\lambda}\ \ \text{and}\ \ Var(\mathcal{E}(\widehat{\textbf{p}}_n)\leq O\left( \frac{(\log n)^2}{n} \right)$$ provided that the probability distribution $(p_j)_{j\in J}$ satisfies $p_j=C_\lambda j^{-\lambda},\,$ where $\lambda>1$.\\

\bigskip \noindent 
Under distributions $p_j=Cj^{-\lambda}$, a necessary condition for$$
\sqrt{n}(\mathcal{E}(\widehat{\textbf{p}}_n)-\mathcal{E}(\textbf{p}))$$ to hold
asymptotic normality is $\lambda\geq 2$.
 \\
 
 \bigskip 
 \noindent \cite{jaya} focused on the number of samples needed to estimate the $\alpha-$Reyni entropy. \\
 
 % voir fich reyni-intro-estimation\\

 \bigskip \noindent However, to our knowledge, no results
regarding the almost sure consistency and the asymptotic normality of the most of entropies, are known.  
\subsection{Main contribution}
 $ $ \\

\noindent Most texts on entropy estimation deal with Shannon entropy estimation whereas we deal with estimation of the most common entropies including Shannon, Tsallis, Reyni, Landsberg-Vedral, Abel entropies, etc by deriving their almost sure convergence and central limit Theorems.\\

\bigskip \noindent Our method consist in getting first general laws for an arbitrary summation of the form 

\begin{equation} \label{defJ}
J(\textbf{p})=\sum_{j\in J} \phi(p_j),
\end{equation}
\noindent where $\phi : (0,1)\rightarrow \mathbb{R}$ is a twice continuously differentiable function.\\

\noindent The results on the summation $J(\textbf{p})$, which is also known under the name of \textit{$\phi$-entropy summation},  will lead to results of entropies already mentioned above.
%\subsection{Paper outline}
\subsection{Overview of the paper}
$ $ \\

\noindent The rest of the paper is organized as follows. In \textsc{S}ection \ref{empifonct}, we define estimators $p_n^{c_j}$ of the \textit{p.m.f} $p_j$ and construct the plug-in estimators of the $\phi-$ entropy summation $J(\textbf{p})=\sum_{j\in J}\phi(p_j)$, where $\phi$ is a twice continously differentiable function, from an i.i.d. sample of size $n$ and according to $\textbf{p}$. We end this section by giving our full results for the summation $J(\textbf{p})$.\\

\noindent In \textsc{S}ection \ref{particular-entropy}, we will particularize the results for specific entropies we already described. \textsc{S}ection \ref{proofentr} provides  the proofs and in \textsc{S}ection \ref{simulation} we present some simulations confirming our results. Finally,  
in \textsc{S}ection \ref{conclusion}, we conclude.\\

\section{ $\phi-$Entropy summation}\label{empifonct}
\subsection{Notations and main results}\label{notation_mainresult}
$ $ \\

\noindent Let $X$ be a random variable defined on the probability  space $(\Omega,\mathcal{A},\mathbb{P})$ and taking values $\mathcal{X}=\{c_1,c_2,\cdots,c_r\}$ with \textit{p.m.f.} $\mathbf{p}=(p_j)_{1\leq j\leq r}$ i.e, 
$$p_j=\mathbb{P}(X=c_j),\ \ \forall j\in J=\{1,2,\cdots,r\}.$$

\bigskip \noindent In general, the full probability distribution $\mathbf{p}=(p_j)_{1\leq j\leq r}$  is not known
and, in particular, in many situations only sets from which to infer entropies are
available. \\

\noindent For example, it could be of interest to determine the  entropies of a given
$DNA$ sequence. In such a case, one could estimate the probability of each element $c_i$ to occur, $p_i$.\\

\noindent Let $X_1,\cdots,X_n$ be $n$ i.i.d. random variables according to $\textbf{p}$. For a given $j\in J$, define the easiest and most objective estimator of $p_j$, based on the i.i.d sample $ X_1,\cdots,X_n,$ by 
 \begin{eqnarray}\label{pn}
 \widehat{p}_n^{c_j}&=&\frac{1}{n}\sum_{i=1}^n1_{c_j}(X_i)%=\frac{1}{n}\#\{i\in \{1,\cdots,n\},\ \ X_i=c_j\}%=\mathbb{P}_{n,X}(1_{c_j})
\end{eqnarray}where 
 $1_{c_j}(X_i)=\begin{cases}
 1\ \ \text{if}\ \ X_i=c_j\\
 0\ \ \text{otherwise}
 \end{cases} $ for any $j\in J$.\\

\noindent For a given $j\in J$, this empirical estimator $\widehat{p}_{n}^j$ of $ p_{j} $ is strongly consistent and asymptotically normal. Precisely, when $n$ tends to infinity,
  \begin{eqnarray} \label{pxn}
  && \widehat{p}_{n}^j-  p_{j} \stackrel{a.s.}{\longrightarrow} 0\\ 
 &&\label{pnj}  \sqrt{n}(\widehat{p}_{n}^j-  p_{j})  \stackrel{\mathcal{D}}{\rightsquigarrow}Z_{p_j}.
\end{eqnarray}  where $Z_{p_j}\stackrel{d }{\sim}\mathcal{N}(0, p_{j} (1- p_{j} ))$

\noindent We denote by $\stackrel{a.s.}{ \longrightarrow}$ the \textit{almost sure convergence} and  $\stackrel{\mathcal{D}}{ \rightsquigarrow}
$ the \textit{convergence in distribution}. The notation $\stackrel{d }{\sim}$ denote the \textit{equality in distribution}. \\

\noindent These asymptotic properties derive from the law of large
numbers and central limit theorem.
\\

\bigskip \noindent The entropy of $\textbf{p}$ can be 
approximated by simply replacing the probabilities $p_j$ by $\widehat{p}_n^{c_j}$ in the entropy summation. For
example, the Shannon entropy $\mathcal{E}_{Sh}(\textbf{p})$ can be estimated by its counter part plug-in
$$\mathcal{E}_{Sh}(\widehat{\textbf{p}}_{n})= -\sum_{j=1}^r \widehat{p}_n^{c_j}\ln (\widehat{p}_n^{c_j})$$

\subsection{$\phi-$entropy summation}

\label{phi_divergence}
\begin{definition}Let $\phi : (0,1)\rightarrow \mathbb{R}$ a twice continously differentiable function.
 The $\phi$-entropy summation of the probability distribution $\textbf{p}=(p_j)_{j\in J}$ is given by
 \begin{equation}
J(\textbf{p})=\sum_{j\in J}\phi(p_j). 
 \end{equation}
 \end{definition}
\noindent  The results on the summation $J(\textbf{p})$ will lead to those on the particular cases of the \textit{Shannon}, \textit{Reyni},  \textit{Tsallis}, \textit{Landsberg-Vedral, Abe}, \textit{Varma} and \textit{$\kappa$-entropies}. 

%
%\begin{table}
%\begin{tabular}{ c c}
%&sd\\
%\hline
%\end{tabular}
%\end{table}
 
\bigskip \noindent  Based on \eqref{pn}, we will use the following $\phi$-entropy summation. 
\begin{eqnarray*}
J(\widehat{\textbf{p}}_n)=\sum_{j\in J}\phi(\widehat{p}_n^{c_j}).
\end{eqnarray*}

\subsection{Statement of the main result}\label{mainresul}
 \noindent It concerns the almost sure efficiency and the asymptotic normality of the summation $\phi$-entropy $J(\widehat{\textbf{p}}_n)$.\\

\noindent Denote  \begin{eqnarray*}
A_J(\textbf{p})&=&\sum_{j\in J}|\phi'( p_{j} )|\\
\text{and}\ \ \sigma^2(\textbf{p})&=& \sum_{j\in J} p_j(1-p_j)( \phi'( p_{j} ))^2-2\sum_{(i,j)\in J^2,i\neq j}(p_ip_j)^{3/2}\phi'( p_{i} )\phi'( p_{j} )\end{eqnarray*}

\begin{theorem} \label{theoj} Let $\textbf{p}=(p_j)_{j\in J}$ a probability distribution and  
$\widehat{\textbf{p}}_n=(\widehat{p}_n^{c_j} )_{j\in J}$ be generated by i.i.d. sample $X_1,\cdots,X_n$ copies of a random variable $X$ according to $\textbf{p}$ and \eqref{BD} be satisfied. 
Then the following asymptotic results hold

\begin{eqnarray}\label{theojas}
 \limsup_{n\rightarrow+\infty}\frac{ \left\vert J(\widehat{\textbf{p}}_{n})-J(\textbf{p})\right \vert }{ a_n }\leq A_J(\textbf{p}),\ \ \text{a.s.},
 \end{eqnarray}

 \begin{equation}\label{theojan}
  \sqrt{n}( J(\widehat{\textbf{p}}_{n})-J(\textbf{p}))\stackrel{\mathcal{D}}{ \rightsquigarrow} \mathcal{N}\left(0,\sigma_J^2(\textbf{p}) \right),\ \ \text{as}\ \ n\rightarrow+\infty.
 \end{equation}
\end{theorem}

\section{ Entropies asymptotic limit law}\label{particular-entropy}

\noindent (\textbf{A}-) Asymptotic behavior of $\mathcal{S}_\alpha(\widehat{\textbf{p}}_n)$.
$\,$\\

\noindent For $\alpha\in (0,1)\cup (1,+\infty)$, denote
 
\begin{eqnarray*}
A_{\mathcal{S}_\alpha}(\textbf{p})&=&\alpha\sum_{j\in J} p_{j}^{\alpha-1}\\
\text{and}\ \ \sigma_{\mathcal{S}_\alpha}^2(\textbf{p})&=&\alpha^2\left( \sum_{j\in J} (1- p_{j} ) p_{j}^{2\alpha-1}-2\sum_{(i,j)\in J^2,i\neq j}(p_ip_j)^{\alpha+1/2}\right).
\end{eqnarray*} 

\begin{corollary}\label{corSalp}
Under the same assumptions as in Theorem \ref{theoj} and for $\alpha\in (0,1)\cup (1,+\infty)$, the following hold
\begin{eqnarray}\label{ialphps}
&& \limsup_{n\rightarrow+\infty}\frac{|\mathcal{S}_\alpha( \widehat{\textbf{p}}_n )-\mathcal{S}_\alpha(\textbf{p})|}{ a_n }\leq A_{\mathcal{S}_\alpha}(\textbf{p}),\ \ \text{a.s.}\\
&&\sqrt{n}(\mathcal{S}_\alpha( \widehat{\textbf{p}}_n )-\mathcal{S}_\alpha(\textbf{p}))\stackrel{\mathcal{D} }{\rightsquigarrow}\mathcal{N}(0,\sigma_{\mathcal{S}_\alpha}^2(\textbf{p})),\ \ \text{as}\ \ n\rightarrow+\infty.\label{ialphno}
\end{eqnarray}
\end{corollary}

\bigskip

\bigskip \noindent (\textbf{B)-} 
 Asymptotic behavior of \textit{Shannon entropy estimator.}
$ $ \\

\bigskip \noindent Let 
\begin{eqnarray*}
A_{Sh}(\textbf{p})&=&\sum_{j\in J}\left\vert 1+ \ln  ( p_{j} )\right \vert\\
\text{and}\ \ \ 
\sigma_{Sh}^2(\textbf{p})&=&\sum_{j\in J} p_j(1- p_{j} )(1+ \ln  ( p_{j} ))^2-2\sum_{(i,j)\in J^2,i\neq j}(p_ip_j)^{3/2}(1+\ln (p_i))(1+\ln (p_j)).\end{eqnarray*}

\begin{corollary}
 \label{theshan}
Under the same assumptions as in Theorem \ref{theoj}, the following hold
\begin{eqnarray}
\limsup_{n\rightarrow +\infty}\frac{| \mathcal{E}_{Sh}(\widehat{p}_{n})-\mathcal{E}_{Sh}(\textbf{p})|}{ a_n }\leq A_{Sh}(\textbf{p}),\ \ \text{a.s.}.
\end{eqnarray}

\begin{eqnarray}
\sqrt{n}\left(\mathcal{E}_{Sh}( \widehat{\textbf{p}}_n )-\mathcal{E}_{Sh}(\textbf{p}) \right)\stackrel{\mathcal{D}}{ \rightsquigarrow} \mathcal{N}(0,\sigma_{Sh}^2(\textbf{p})),\ \ \text{as}\ \ n\rightarrow +\infty.
\end{eqnarray}

\end{corollary}

\bigskip
\noindent 
( \textbf{C-)} Asymptotic behavior of the \textit{Reyni entropy  estimator}.$ $ \\ 
 
 \bigskip \noindent The treatment of the asymptotic behavior of the Renyi-$\alpha$ entropies estimator and of the $\alpha,\beta-$Varma entropy estimator is obtained by the application of the delta method.\\

\noindent For $\alpha\in (0,1)\cup (1,+\infty)$, denote 
\begin{eqnarray*}
A_{R,\alpha}(\textbf{p})&=& \frac{\alpha}{\left\vert\alpha-1\right\vert \mathcal{S}_\alpha(\textbf{p}) } \sum_{j\in J}  p_{j} ^{\alpha-1}\\
\text{and}\ \  \ \sigma _{R,\alpha}^{2}(\textbf{p})&=&\left(\frac{\alpha}{(\alpha-1)\mathcal{S}_\alpha(\textbf{p})}\right)^2\left( \sum_{j\in J} (1- p_{j} ) p_{j}^{2\alpha-1}-2\sum_{(i,j)\in J^2,i\neq j} (p_ip_j)^{\alpha+1/2}\right).
\end{eqnarray*}

\begin{corollary}\label{corRen} Under the same assumptions as in Theorem \ref{theoj} and for any $\alpha\in (0,1)\cup (1,+\infty)$, the following hold
\begin{eqnarray}\label{corRenas}
&& \limsup_{n\rightarrow +\infty}\frac{|\mathcal{E}_{R,\alpha}( \widehat{\textbf{p}}_n )-\mathcal{E}_{R,\alpha}(\textbf{p})|}{ a_n }\leq A_{R,\alpha}(\textbf{p}),\ \ \text{a.s}.
\end{eqnarray}

\begin{eqnarray}\label{corRenan}
\sqrt{n}\left( \mathcal{E}_{R,\alpha}( \widehat{\textbf{p}}_n )-\mathcal{E}_{R,\alpha}(\textbf{p})\right)\stackrel{\mathcal{D}}{ \rightsquigarrow} & 
\mathcal{N}\left( 0,\sigma _{\mathcal{R}\alpha}^{2}(\textbf{p})\right)\text{ as } n\rightarrow + \infty.
\end{eqnarray}
\end{corollary}

\bigskip
\noindent (\textbf{D}-) Asymptotic behavior of the \textit{Tsallis entropy estimator}.$ $ \\ 

\bigskip \noindent For $\alpha\in (0,1)\cup (1,+\infty)$, denote 

\begin{eqnarray*}
A_{T,\alpha}(\textbf{p})&=& \frac{\alpha}{| \alpha-1|}\sum_{j\in J}
 p_{j} ^{\alpha-1}\\
 \text{and}\ \ \sigma_{T,\alpha}^2(\textbf{p})&=&\left(\frac{\alpha}{\alpha-1}\right)^2 \left( \sum_{j\in J} (1- p_{j} ) p_{j}^{2\alpha-1}-2\sum_{(i,j)\in J^2,i\neq j}(p_ip_j)^{\alpha+1/2}\right).
\end{eqnarray*}

\begin{corollary}\label{corTs}
Under the same assumptions as in Theorem \ref{theoj} and for $\alpha\in (0,1)\cup (1,+\infty)$, the following hold
\label{corts1}
\begin{eqnarray} \limsup_{n\rightarrow +\infty}\frac{|\mathcal{E}_{T,\alpha}( \widehat{\textbf{p}}_n )-\mathcal{E}_{T,\alpha}(\textbf{p})|}{ a_n }\leq A_{T,\alpha}(\textbf{p})\ \ \text{a.s}.
\end{eqnarray}\text

\begin{eqnarray}\label{corts2}
\sqrt{n}\left(\mathcal{E}_{T,\alpha}( \widehat{\textbf{p}}_n )-\mathcal{E}_{T,\alpha}(\textbf{p}) \right)\stackrel{\mathcal{D}}{ \rightsquigarrow} \mathcal{N}(0,\sigma_{T,\alpha}^2(\textbf{p}))\ \ \text{as}\ \ n\rightarrow+\infty.
\end{eqnarray}
\end{corollary}

\bigskip 

\noindent (\textbf{E-}) Asymptotic behavior of the \textit{Landsberg-Vedral entropy} estimator.$\,$\\

\bigskip \noindent For $\alpha\in (0,1)\cup (1,+\infty)$, denote 

\begin{eqnarray*}
A_{L.V,\alpha}(\textbf{p})&=& \frac{\alpha}{| \alpha-1|\mathcal{S}_\alpha(\textbf{p}) }\sum_{j\in J}
 p_{j} ^{\alpha-1} \\
 \text{ and}\ \ \
   \sigma_{L.V,\alpha}^2(\textbf{p})&=&\left(\frac{\alpha}{(\alpha-1)\mathcal{S}_\alpha(\textbf{p})}\right)^2 \left( \sum_{j\in J} (1- p_{j} ) p_{j}^{2\alpha-1}-2\sum_{(i,j)\in J^2,i\neq j}(p_ip_j)^{\alpha+1/2}\right).
\end{eqnarray*}

 \begin{corollary}\label{corlandval}
Under the same assumptions as in Theorem \ref{theoj} and for $\alpha\in (0,1)\cup (1,+\infty)$, the following hold \begin{eqnarray}
&&\limsup_{n\rightarrow+\infty}\frac{|\mathcal{E}_{L.V,\alpha}( \widehat{\textbf{p}}_n )-\mathcal{E}_{L.V,\alpha}( \textbf{p})|}{ a_n }\leq A_{L.V,\alpha}(\textbf{p}),\ \ \text{a.s.}\\
&& \sqrt{n}\left( \mathcal{E}_{L.V,\alpha}( \widehat{\textbf{p}}_n )-\mathcal{E}_{L.V,\alpha}(\textbf{p})\right)\stackrel{\mathcal{D}}{\rightsquigarrow} \mathcal{N}\left(0, \sigma_{L.V,\alpha}^2(\textbf{p})\right)\ \ \text{as}\ \ n\rightarrow+\infty.
\end{eqnarray}
\end{corollary}

\bigskip 
\noindent (\textbf{F}-)  Asymptotic behavior of $\alpha-$\textit{Abel entropy estimator}.$ $ \\ 

%
%\noindent Denote  \begin{eqnarray*}
%A_J(\textbf{p})&=&\sum_{j\in J}|\phi'( p_{j} )|\\
%\text{and}\ \ \sigma^2(\textbf{p})&=& \sum_{j\in J} p_j(1-p_j)( \phi'( p_{j} ))^2-2\sum_{(i,j)\in J^2,i\neq j}(p_ip_j)^{3/2}\phi'( p_{i} )\phi'( p_{j} )\end{eqnarray*}
% $$f(x)=-\frac{\alpha}{\alpha^2-1}(x^\alpha-x^{\alpha^{-1}})\ \ \ f'(x)=-\frac{1}{\alpha^2-1}(\alpha^2 x^{\alpha-1}- x^{(1/\alpha)-1})$$
  \noindent For $\alpha\in (0,1)\cup (1,+\infty)$, denote 
\begin{eqnarray*}
A_{\mathcal{A}b,\alpha}(\textbf{p})&=&\frac{1}{|\alpha^2-1|} \sum_{j\in J}\left\vert \alpha^2  p_{j} ^{\alpha-1}-   p_{j} ^{(1/\alpha)-1}\right\vert \ \ \\
\text{and}\ \ \sigma _{\mathcal{A}b,\alpha}^{2}(\textbf{p})&=&\frac{1}{(\alpha^2-1)^2}\biggr( \sum_{j\in J} (1- p_{j} )\left( \alpha^2  p_{j} ^{\alpha-1/2}-   p_{j} ^{(1/\alpha)-1/2}\right)^2\\
&&\ \ \ \ -2\sum_{(i,j)\in J^2,i\neq j}\left[\alpha^2  p_{i} ^{\alpha+1/2}-   p_{i} ^{(1/\alpha)+1/2}\right]\left[\alpha^2  p_{j} ^{\alpha+1/2}-   p_{j} ^{(1/\alpha)+1/2}\right]\biggr).
\end{eqnarray*}

\begin{corollary}
\label{corabel}
Under the same assumptions as in Theorem \ref{theoj} and for any $\alpha\in (0,1)\cup (1,+\infty)$, the following hold \begin{eqnarray}
&&\limsup_{n\rightarrow+\infty}\frac{|\mathcal{E}_{Ab,\alpha}( \widehat{\textbf{p}}_n )-\mathcal{E}_{Ab,\alpha}( \textbf{p })|}{ a_n }\leq A_{\mathcal{A}b,\alpha}(\textbf{p})\\
&& \sqrt{n}\left( \mathcal{E}_{Ab,\alpha}( \widehat{\textbf{p}}_n )-\mathcal{E}_{Ab,\alpha}( \textbf{p })\right)\stackrel{\mathcal{D}}{\rightsquigarrow} \mathcal{N}(0,\sigma _{\mathcal{A}b,\alpha}^{2}(\textbf{p})),\ \ \text{as}\ \ n\rightarrow +\infty.
\end{eqnarray}
\end{corollary}

\bigskip

\noindent (\textbf{G-}) Asymptotic \textit{behavior of 
$\kappa-$entropy}
$\,$
\\

 \noindent  For $\kappa\in (0,1),$ denote 
\begin{eqnarray*}
A_{\kappa}(\textbf{p})&=&\frac{1}{2\kappa }\sum_{j\in J}\left\vert (1-\kappa) p_{j} ^{-\kappa}-(1+\kappa) p_{j} ^{\kappa}\right\vert\\
\text{and}\ \ \sigma_{\kappa}^{2}(\textbf{p})&=&\frac{1}{4\kappa^2 }\biggr(\sum_{j\in J}  (1- p_{j} ) \left( (1-\kappa) p_{j} ^{-\kappa+1/2}-(1+\kappa) p_{j} ^{\kappa+1/2}\right)^2\\
&& -2\sum_{(i,j)\in J^2,i\neq j}\left[(1-\kappa) p_{i} ^{-\kappa+3/2}-(1+\kappa) p_{i} ^{\kappa+3/2}\right]\,\left[(1-\kappa) p_{j} ^{-\kappa+3/2}-(1+\kappa) p_{j} ^{\kappa+3/2}\right]\biggr).
\end{eqnarray*}

\begin{corollary}\label{corkappa}
Under the same assumptions as in Theorem \ref{theoj} and for any $\kappa\in (0,1)$, the following hold
\begin{eqnarray*}
&&\limsup_{n\rightarrow+\infty}\frac{ |\mathcal{E}_{\kappa}( \widehat{\textbf{p}}_n )- \mathcal{E}_{\kappa}(\textbf{p})|}{ a_n }\leq A_{\kappa}(\textbf{p}),\ \ \text{a.s}\\
&& \sqrt{n}(\mathcal{E}_{\kappa}( \widehat{\textbf{p}}_n )- \mathcal{E}_{\kappa}(\textbf{p}))\stackrel{\mathcal{D} }{\rightsquigarrow}
\mathcal{N}(0,\sigma_{\kappa}^{2}(\textbf{p})),\ \ \text{as}\ \ n \rightarrow +\infty.\end{eqnarray*}
\end{corollary}
\bigskip

\noindent ( \textbf{H }-) 
Asymptotic behavior of Varma's entropy of order $\alpha$ and type $\beta$.\\

 \noindent For $\beta-1<\alpha<\beta,\ \ \beta\geq 1$ denote \begin{eqnarray*}
A_{V,\alpha,\beta}(\textbf{p})&=&\frac{  \alpha+\beta-1}{S_{\alpha+\beta-1}  } \sum_{j\in J} p_{j}^{\alpha+\beta-2}\\
\text{and}\ \ \sigma_{V,\alpha,\beta}^2(\textbf{p})
&=&\left(
\frac{\alpha+\beta-1}{(\beta-\alpha)\mathcal{S}_{\alpha+\beta-1}(\textbf{p})}\right)^2
\left( \sum_{j\in J} (1- p_{j} ) p_{j}^{2\alpha+2\beta-3}-2\sum_{(i,j)\in J^2,i\neq j}(p_ip_j)^{\alpha+\beta+1/2}
\right).
\end{eqnarray*}

\begin{corollary}\label{varma}
Under the same assumptions as in Theorem \ref{theoj} and for $\beta-1<\alpha<\beta,\ \ \beta\geq 1$, the following hold
\begin{eqnarray}
&& \label{varmas}\limsup_{n\rightarrow+\infty}\frac{|\mathcal{E}_{V,\alpha,\beta} ( \widehat{\textbf{p}}_n )-\mathcal{E}_{V,\alpha,\beta}|}{ a_n }\leq A_{V,\alpha,\beta}(\textbf{p})
% \frac{1}{|1-\alpha|}A_{\mathcal{S}_\alpha}
,\ \ \text{a.s.}\\
&&\label{smna} \sqrt{n}( \mathcal{E}_{V,\alpha,\beta} ( \widehat{\textbf{p}}_n )-\mathcal{E}_{V,\alpha,\beta} (\textbf{p}))\stackrel{\mathcal{D} }{\rightsquigarrow}\mathcal{N}(0,\sigma_{V,\alpha,\beta}^2(\textbf{p})),\ \ \text{as}\ \ n\rightarrow+\infty.
\end{eqnarray}
\end{corollary}

\section{The proofs}\label{proofentr}
\noindent Before we state our main results we introduce the following  notations. For a fixed $j\in J$, denote
\begin{eqnarray*}
 \Delta_{p_n}^{c_j} =\widehat{p}_n^{c_j}-  p_{j} ,\ \ 
\delta_n(p_j)=\sqrt{n/p_j}\Delta_{p_n}^{c_j} , 
\end{eqnarray*}
and $
a_n =\sup_{j\in J}|  \Delta_{p_n}^{c_j}  |.$\\

\noindent We recall that,
since for a fixed $j\in D,$ $n\widehat{p}_n^{c_j}$ has a binomial distribution with parameters $n$  and success probability $p_j$, we have 
 \begin{equation*}
 \mathbb{E}\left[ \widehat{p}_n^{c_j}\right]=p_j\ \ \text{and}\ \ \mathbb{V}(\widehat{p}_n^{c_j})=\frac{p_j(1-p_j)}{n}.
\end{equation*}

\noindent And finally, by the asymptotic Gaussian limit of the multinomial law (see for example \cite{ips-wcia-ang}, Chapter 1, Section 4), we have
\begin{eqnarray}
\label{pnj}&& \biggr( \delta_n(p_j), \ j\in J\biggr)
\stackrel{\mathcal{D}}{\rightsquigarrow }Z(\textbf{p}) \stackrel{d }{\sim}\mathcal{N}(0,\Sigma_\textbf{p}),\ \ \text{as}\ \ n\rightarrow +\infty,
\end{eqnarray}where $Z(\textbf{p})= (Z_{p_j},j\in J)^t
$ \
is a centered Gaussian random vector of dimension $\#(J)$ having the following elements :
\begin{eqnarray}\label{vars}
&&\left(\Sigma_\textbf{p}\right)_{(i,j)}=(1-p_j)\delta_{ij}-\sqrt{ p_ip_j} (1-\delta_{ij}), \ \ (i,j) \in J^2,
\end{eqnarray}where  $\delta_{ij}=\begin{cases}
1\ \ \text{for}\ \ i=j\\
0\ \ \text{for}\ \ i\neq j
\end{cases}.$

\subsection{Proof of Theorem \ref{theoj}}

For a fixed $j\in J$, we have
\begin{eqnarray}\label{eq1}
\notag \phi( \widehat{p}_n^{c_j} )&=&\phi( p_{j} + \Delta_{p_n}^{c_j} )\\
&=&\phi( p_{j} )+ \Delta_{p_n}^{c_j} \phi'( p_{j} +\theta_{1}(j) \Delta_{p_n}^{c_j} ),
\end{eqnarray} by the mean value Theorem applied to the function $\phi$ and where $\theta_{1}(j)\in (0,1).$ \\
Apply again the mean value Theorem to the derivative of the function $\phi'$
\begin{eqnarray*}
\phi'( p_{j} +\theta_{1}(j) \Delta_{p_n}^{c_j} )&=&\phi'( p_{j} )+\theta_{1}(j) \Delta_{p_n}^{c_j} \phi"( p_{j} +\theta_{2}(j) \Delta_{p_n}^{c_j} ),
\end{eqnarray*}where $\theta_{2}(j)\in (0,1).$  We can write \eqref{eq1} as \begin{eqnarray*}
 \phi( \widehat{p}_n^{c_j} )&=&\phi( p_{j} )+  \Delta_{p_n}^{c_j} \phi'( p_{j} )+\theta_{1}(j)( \Delta_{p_n}^{c_j} )^2\phi"( p_{j} +\theta_{2}(j) \Delta_{p_n}^{c_j} )
 \end{eqnarray*}Now we have, by summation over $j\in J$
 \begin{eqnarray}\label{eq2}
J( \widehat{\textbf{p}}_n )-J(\textbf{p})&=&\sum_{j\in J} \Delta_{p_n}^{c_j} \phi'( p_{j} )\\
&& \ \ +\ \ \notag \sum_{j\in J}\theta_{1}(j)( \Delta_{p_n}^{c_j} )^2\phi"( p_{j} +\theta_{2}(j) \Delta_{p_n}^{c_j} )
 \end{eqnarray}Hence 
 \begin{eqnarray*}
 \left\vert J( \widehat{\textbf{p}}_n )-J(\textbf{p})\right \vert &\leq & a_n \sum_{j\in J}|\phi'( p_{j} )|+ a_n ^2\sum_{j\in J}|\phi"( p_{j} +\theta_{2}(j) \Delta_{p_n}^{c_j} )|,
 \end{eqnarray*} Therefore
 \begin{eqnarray*}
 \limsup_{n\rightarrow+\infty}\frac{ \left\vert J( \widehat{\textbf{p}}_n )-J(\textbf{p})\right \vert }{ a_n }\leq A_J(\textbf{p}),\ \ \text{a.s.},
 \end{eqnarray*} since $ a_n 
  \stackrel{a.s.}{\rightarrow}0$  as $n\rightarrow +\infty$ and $$\sum_{j\in J}|\phi"( p_{j} +\theta_{2}(j) \Delta_{p_n}^{c_j} )|\rightarrow \sum_{j\in J}|\phi"( p_{j} )| <\infty,\ \ \text{as}\ \ n \rightarrow +\infty.$$ 
 
\noindent This prove \eqref{theojas}.\\

\noindent Let prove \eqref{theojan}. By going back to \eqref{eq2}, we get \begin{eqnarray*}
\sqrt{n}(J( \widehat{\textbf{p}}_n )-J(\textbf{p}))&=&\sum_{j\in J}\sqrt{p_j} \delta_n(p_j)  \phi'( p_{j} )+\sqrt{n}R_{n},
\end{eqnarray*}where $$R_{n}=\sum_{j\in J}\theta_{1}(j)( \Delta_{p_n}^{c_j} )^2\phi"( p_{j} +\theta_{2}(j) \Delta_{p_n}^{c_j} ).$$

\noindent Using Formula \eqref{pnj} above, we get
 \begin{equation}
\label{pnj2}\sum_{j\in J} \sqrt{p_j} \delta_n(p_j) \phi'( p_{j} )\stackrel{\mathcal{D}}{ \rightsquigarrow} \sum_{j\in J}\  \phi'( p_{j} )\sqrt{p_j}Z_{p_j}
,\ \ \text{as}\ \ n\rightarrow+\infty,
\end{equation}

\noindent which follows a centered normal law of variance $\sigma_J^2( \textbf{p} )$ since 
\begin{eqnarray*}
\text{Var}\left( \sum_{j\in J} \phi'( p_{j} )\sqrt{p_j} Z_{p_j}\right)&=&\sum_{j\in J} \text{Var}\left( \phi'( p_{j} )\sqrt{p_j} Z_{p_j}\right)+2\sum_{j\in J} \text{Cov}\left( \phi'( p_{i} )\sqrt{p_i} Z_{p_i},\phi'( p_{j} )\sqrt{p_j} Z_{p_j}\right)\\
&=& \sum_{j\in J} p_j(1-p_j)( \phi'( p_{j} ))^2-2\sum_{(i,j)\in J^2,i\neq j}p_ip_j\sqrt{p_ip_j}\phi'( p_{i} )\phi'( p_{j} ).
\end{eqnarray*}

\bigskip \noindent The proof will be complete if we show that $\sqrt{n}R_{n}$ converges to zero in probability. \\

\noindent We have  
\begin{equation} \label{r1n}
\left\vert \sqrt{n}R_{n}\right\vert \leq \sqrt{n}a_{n}^{2} \sum_{j\in J}\phi"(p_{j}+\theta _{2,j}\Delta_{p_n}^{c_j}).
\end{equation}
\noindent By the Bienaymé-Tchebychev inequality, we have, for any fixed $\epsilon >0$ and for any $j\in J$
\begin{eqnarray*}
\mathbb{P}(\sqrt{n}( \widehat{p}_n^{c_j} - p_{j} )^2\geq \epsilon)=\mathbb{P}\left(| \widehat{p}_n^{c_j} - p_{j} |\geq  \frac{\sqrt{\epsilon} }{n^{1/4}}\right)\leq \frac{ p_{j} (1- p_{j} )}{\epsilon n^{1/2}}.
\end{eqnarray*} 
\noindent Hence $\sqrt{n}a_{n}^{2}=o_{\mathbb{P}}(1)$, 
%that is 
%$\sqrt{n}
%R_{n}\stackrel{\mathbb{P}}{\rightarrow} 0\ \  \text{ as } \ \ n\rightarrow +\infty$
 which proves \eqref{theojan}. \\

\noindent All this ends 
%implies $$\sqrt{n}(J( \widehat{\textbf{p}}_n )-J(\textbf{p}))\stackrel{\mathcal{D}}{ \rightsquigarrow} \mathcal{N}\left(0,\sigma^2( \textbf{p} ) \right),\ \ \text{as}\ \ n\rightarrow+\infty.$$
the proof of Theorem \ref{theoj}.

\subsection{Proofs of Corollaries} $ $ \\

\noindent \textbf{A}-) 
\noindent The Proofs of Corollaries \ref{corSalp} and \ref{theshan} are direct adaptations of Theorem \ref{theoj} with respectively $\phi(s)=s^\alpha$ and $\phi(s)=-s \ln  s$.\\

\bigskip \noindent \textbf{B}-) Proof of Corollary \ref{corRen}. For $\alpha\in (0,1)\cup (1,+\infty),$ $\alpha-$ Reyni entropy is expressed through the power sum $\mathcal{S}_\alpha(\textbf{p})=\sum_{j\in J}\phi(p_j)$ with $\phi(s)=s^\alpha$. We have
\begin{equation*}
\mathcal{E}_{R,\alpha}( \widehat{\textbf{p}}_n )-\mathcal{E}_{R,\alpha }(\textbf{p})=%
\frac{1}{\alpha -1}\left(  \ln  \mathcal{S}_\alpha( \widehat{\textbf{p}}_n )- \ln \mathcal{S}_\alpha(\textbf{p})\right),
\end{equation*}%
by using a Taylor
expansion of $ \ln  (1+y)$ it follows that, almost surely, 
\begin{eqnarray*}
 \ln  \mathcal{S}_\alpha( \widehat{\textbf{p}}_n )- \ln \mathcal{S}_\alpha(\textbf{p}) &=& \ln  \left( 1+%
\frac{\mathcal{S}_\alpha( \widehat{\textbf{p}}_n )-\mathcal{S}_{\alpha}(\textbf{p})}{\mathcal{S}_{\alpha}(\textbf{p})}%
\right) \\
&=&\frac{\mathcal{S}_\alpha( \widehat{\textbf{p}}_n )-\mathcal{S}_{\alpha}(\textbf{p})}{\mathcal{S}_{\alpha}(\textbf{p})}+O_{\text{a.s%
}}(a_{n}^{2}).
\end{eqnarray*}
\noindent Finally this, combined with \eqref{ialphps} of Corollary \ref{corSalp}, proves \eqref{corRenas}.\\

\bigskip \noindent Now recall by going back to \eqref{eq2}, we can write
\begin{eqnarray*}
\sqrt{n}( \mathcal{S}_\alpha( \widehat{\textbf{p}}_n )-\mathcal{S}_{\alpha}(\textbf{p}))
%&=&\sqrt{n}\sum_{j\in J}( \widehat{p}_n^{c_j} - p_{j} )\phi'( p_{j} )+o_{\mathbb{P}}(1) \\
&=& \sqrt{n}\sum_{j\in J} \Delta_{p_n}^{c_j} \phi'( p_{j} )+o_{\mathbb{P}}(1)
%=O_{\mathbb{P}}(1)
\end{eqnarray*}
here $\phi'( p_{j} )=\alpha  p_{j} ^{\alpha-1}$. \\
\noindent Hence dividing each member by 
$\sqrt{n} \mathcal{S}_{\alpha}(\textbf{p})$, we get\begin{equation*}
\frac{ \mathcal{S}_\alpha( \widehat{\textbf{p}}_n )}{\mathcal{S}_{\alpha}(\textbf{p})}=1+\frac{\sum_{j\in J} \Delta_{p_n}^{c_j} \phi'( p_{j} )}{\mathcal{S}_{\alpha}(\textbf{p})}+o_{\mathbb{P}}(1).
\end{equation*}

\noindent Now by Taylor expansion of $ \ln  (1+y)$, it follows that, almost surely, 
\begin{eqnarray*}
 \ln  \mathcal{S}_\alpha( \widehat{\textbf{p}}_n )- \ln \mathcal{S}_\alpha(\textbf{p}) &=& \ln  \left( 1+%
\frac{\sum_{j\in J} \Delta_{p_n}^{c_j} \phi'( p_{j} )}{\mathcal{S}_\alpha(\textbf{p})}\right) \\
&=&\frac{\sum_{j\in J} \Delta_{p_n}^{c_j} \phi'( p_{j} )}{\mathcal{S}_\alpha(\textbf{p})}+O_{%
\mathbb{P}}\left( \frac{1}{n}\right)
\end{eqnarray*}

\noindent therefore 
\begin{eqnarray*}
\sqrt{n}\left( \mathcal{E}_{R,\alpha}( \widehat{\textbf{p}}_n )-\mathcal{E}_{R,\alpha }(\textbf{p})\right)&=&\frac{1}{\alpha -1}\frac{\sum_{j\in J}\sqrt{n} \Delta_{p_n}^{c_j} \phi'( p_{j} )}{\mathcal{S}_\alpha(\textbf{p})}+o_{\mathbb{P}}(1)\\ &\stackrel{\mathcal{D}}{ \rightsquigarrow} & 
\mathcal{N}\left( 0,\sigma _{R,\alpha}^{2}(\textbf{p})\right)\text{ as } n\rightarrow +\infty,
\end{eqnarray*}using \eqref{pnj2} and
where $$ \sigma _{R,\alpha}^{2}(\textbf{p})=\left(\frac{\alpha}{(\alpha-1)\mathcal{S}_\alpha(\textbf{p})}\right)^2\left( \sum_{j\in J} (1- p_{j} ) p_{j}^{2\alpha-2}-2\sum_{(i,j)\in J^2,i\neq j} (p_ip_j)^{\alpha-1/2}\right).$$

\noindent This proves \eqref{corRenan} and ends the proof of the Corollary \ref{corRen}.  \\

\bigskip \noindent \textbf{C}-) Proof of Corollary \ref{corts1}. Since $\alpha-$Tsallis entropy is related to the power sum $\mathcal{S}_\alpha(\textbf{p})$,  the proof follows directly from Corollary \ref{corSalp}.

\bigskip 
\noindent \textbf{D}-) Proof of Corollary \ref{corlandval}. 
 Since \textit{Landsberg-Vedral} and Tsallis $\alpha-$ entropies are related by 
 
 $$\mathcal{E}_{L.V,\alpha}(\textbf{p})=\frac{\mathcal{E}_{T,\alpha}(\textbf{p})}{ \mathcal{S}_\alpha(\textbf{p})},$$ the proof of this Corollary results directly from the Corollary \ref{corTs}.\\

\noindent \textbf{E-}) The proof of Corollary \ref{varma} is similar to the one of Corollary \ref{corRen} with the power sum $\mathcal{S}_{\alpha+\beta-1}(\textbf{p})=\sum_{j\in J}\phi(p_j)$ with $\phi(s)=s^{\alpha+\beta-1}$.

 \bigskip \noindent \textbf{F-}) Corollaries \ref{corabel} and \ref{corkappa} are, as for Corollaries \ref{corSalp} and \ref{theshan},  adaptations of Theorem \ref{theoj} with this time $\phi(s)=\frac{-1}{\alpha-\alpha^{-1}}(s^\alpha-s^{\alpha^{-1}})$ and $\phi(s)=\frac{1}{2\kappa } (s^{1+\kappa}- s^{1-\kappa})
,\ \ \kappa\in (0,1),$ respectively.

\section{Simulation}\label{simulation}
\noindent \noindent 
To assess the performance of ours estimators, we present a simulation study.  \\

 \noindent  Let  $X$ a random variable defined on a measurable space $(\Omega,\mathcal{A},\mathbb{P})$ and with range $\mathcal{X}=\{1,2,3\}$ with their respective probabilities mass $$
 p_1=0.4,\ \ p_2=0.25,\ \ p_3=0.35.$$  

\noindent We plot the entropies estimators and construct histograms and Q-Q plots to see whether data are normally distributed.\\

\noindent In each figure, the left panel represents
 the plot of entropy estimator, built from sample sizes of $n=100,200,\cdots,30000$, and the true entropy (represented by horizontal black line). The
middle panel shows the histogram of the  data and the red line represents the plot of the theoretical normal distribution calculated
from the same mean and the same standard deviation of the data. The right panel concerns the Q-Q plot of the data which display the observed values against normally 
distributed data (represented by the red line).\\

\begin{figure}[H]
\centering
\includegraphics[scale=0.3]{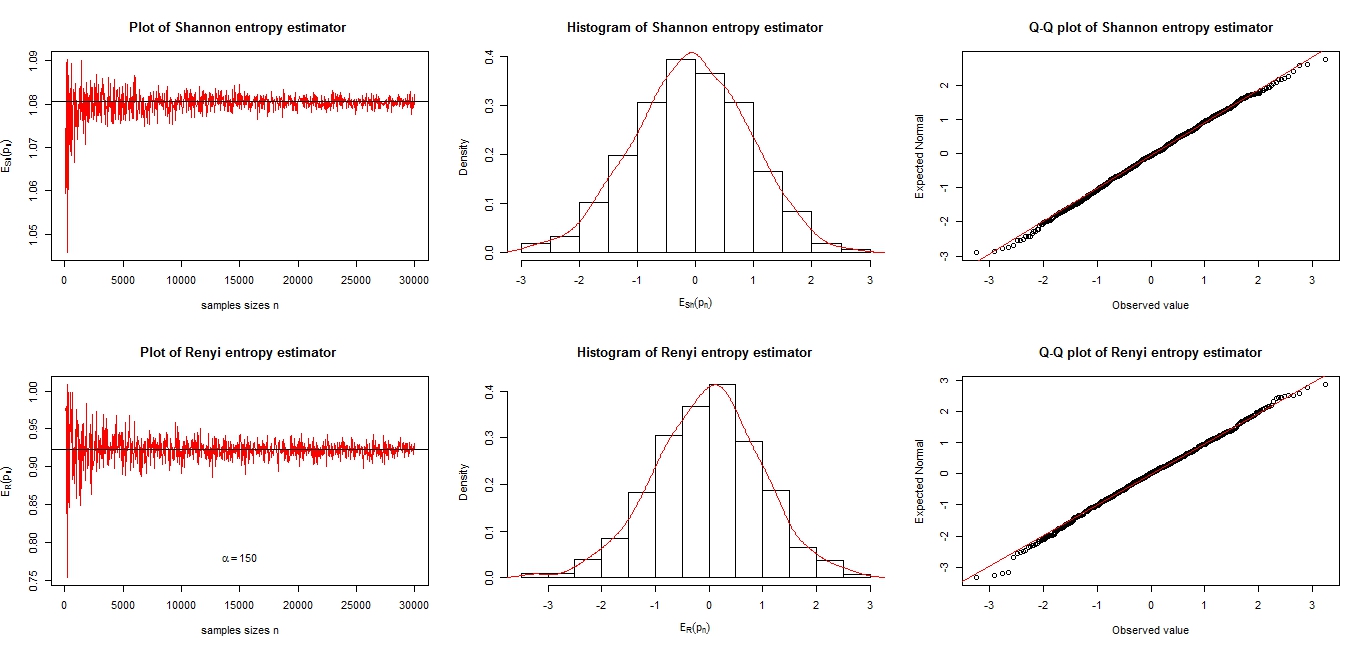} 
 \caption{Plots of \textit{Shannon} and \textit{Renyi} entropies estimators when samples sizes increase, histograms and normal Q-Q plots  versus $\mathcal{N}(0,1)$.
}\label{shanrreny}
\end{figure}

\begin{figure}[H]
\centering
\includegraphics[scale=0.3]{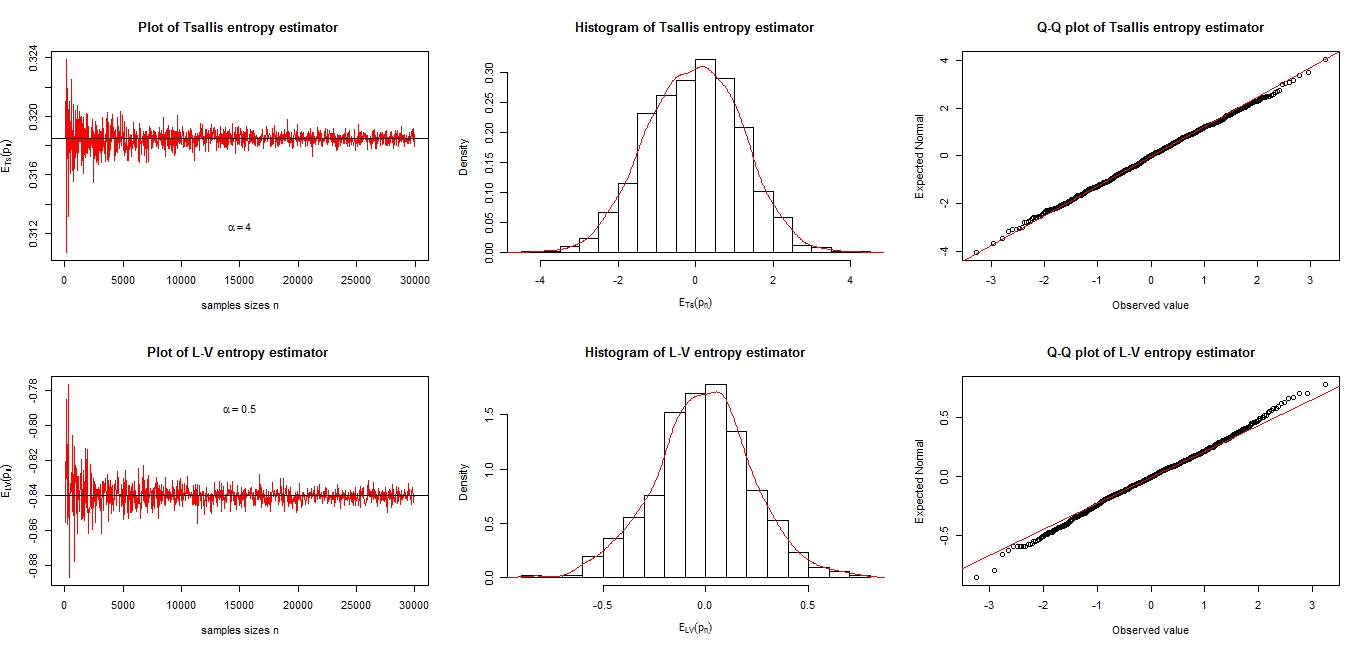} 
 \caption{Plots of \textit{Tsallis} and \textit{Landsberg-Vedral} entropies estimators when samples sizes increase, histograms and normal Q-Q plots  versus $\mathcal{N}(0,1)$.}\label{tsallv}
\end{figure}

\begin{figure}[H]
\centering
\includegraphics[scale=0.3]{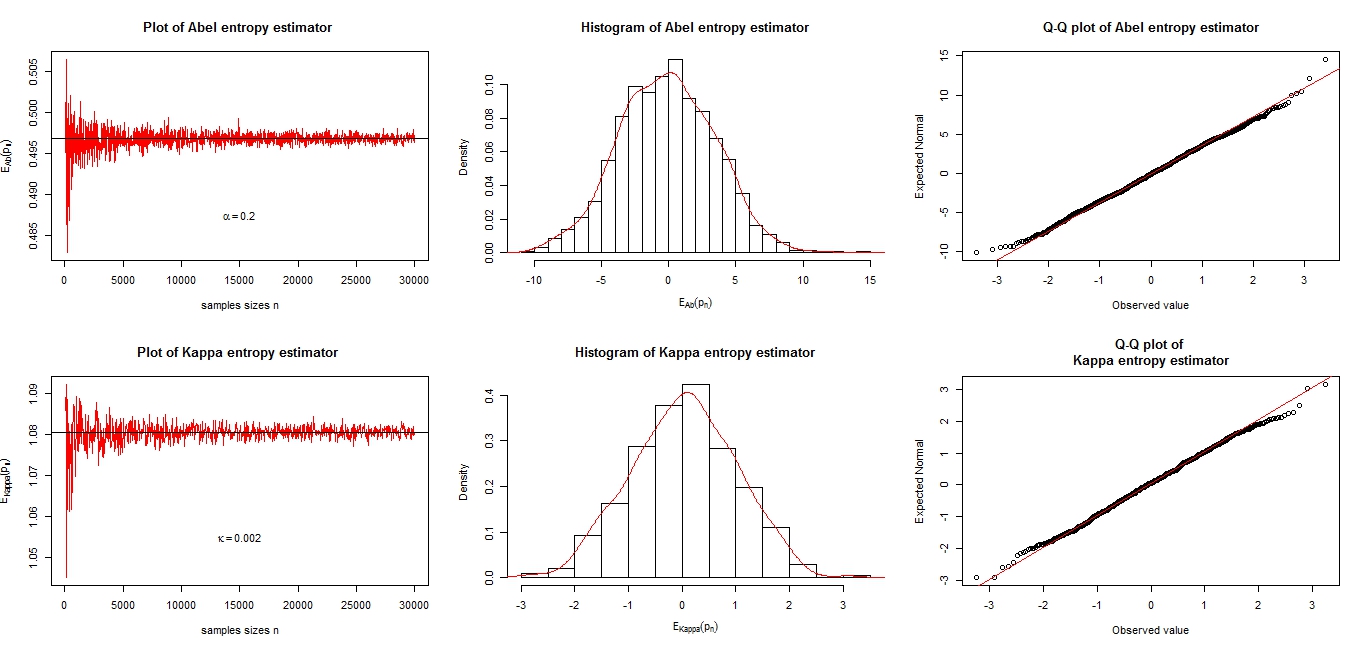} 
 \caption{Plots of \textit{Abel} and \textit{Kappa} entropies estimators when samples sizes increase, histograms and normal Q-Q plots  versus $\mathcal{N}(0,1)$.}\label{AbeKap}
\end{figure}

\begin{figure}[H]
\centering
\includegraphics[scale=0.3]{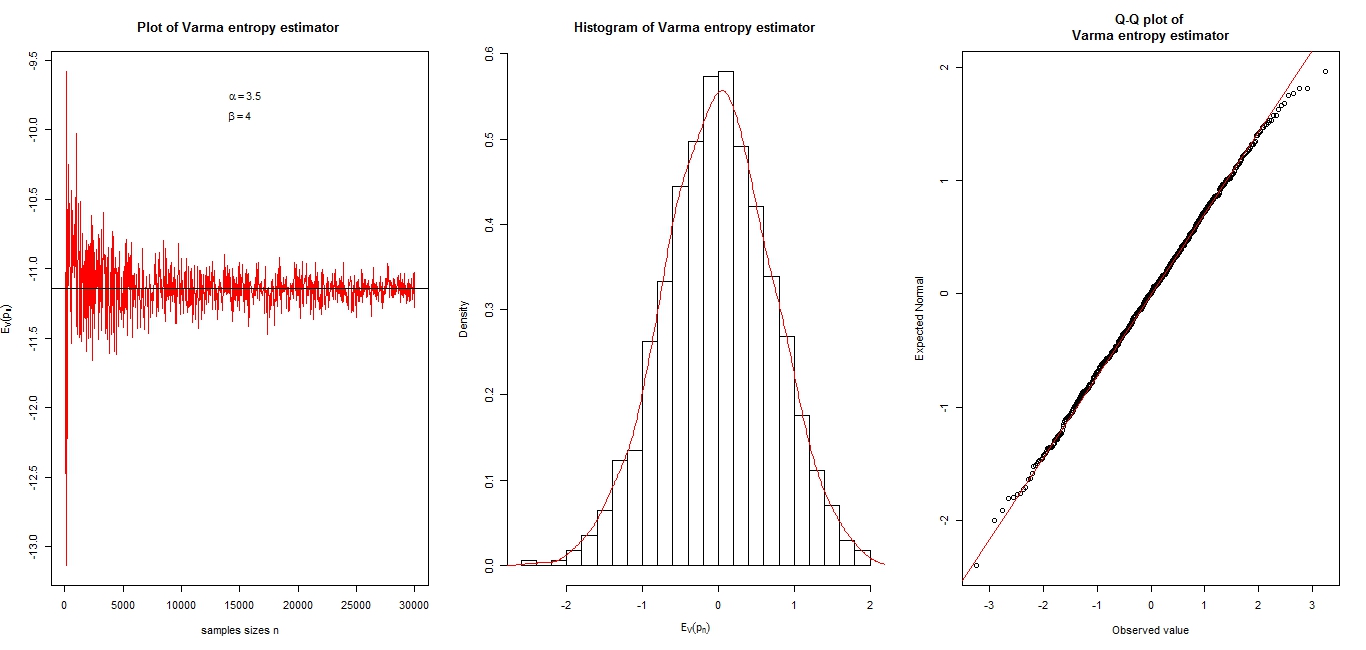}  \caption{Plot of \textit{Varma} entropy estimator when samples sizes increase, histogram and normal Q-Q plots  versus $\mathcal{N}(0,1)$.
}\label{varm}
\end{figure}

\noindent As we can see from \textsc{figures}  \ref{shanrreny} \ref{tsallv}, \ref{AbeKap} and 
\ref{varm}, our entropies estimators are asymptotically normally 
distributed. 

\section{Conclusion}\label{conclusion}

\noindent We have derived a new nonparametric estimator for
 entropies in
the discrete case and on finite sets. We adopted the plug-in method
and we derived almost sure rates of convergence and central limit Theorems for some of the most important entropies in the discrete case. We also
demonstrated their efficiency using a simulation study.

\end{document}